\documentclass[runningheads]{llncs}
\usepackage[T1]{fontenc}
\usepackage{graphicx}
\usepackage{comment}
\usepackage{subcaption}
\usepackage{amsmath}
\usepackage{amssymb, amsfonts, mathrsfs}
\usepackage{hyperref}
\usepackage{xcolor}

\begin{document}
\title{Portfolio Optimization with `Physical' Decision Variables and Non-Linear Performance Metrics: Diversification Challenge and Proposals}
%
%
\author{Isabel Barros Garcia \inst{1} \and
Jérémie Messud \inst{1}}
\authorrunning{I. Barros Garcia and J. Messud}
%
\institute{TotalEnergies, Tour Coupole - 2 place Jean Millier 92078 Paris la Defense cedex, France \\ 
\email{\{isabel.barros-garcia,jeremie.messud\}@totalenergies.com}}
\maketitle              
\begin{abstract}
Portfolio optimization (PO) is a core tool in financial and operational decision-making, typically balancing expected profit and risk. In real-world applications, particularly in the energy sector, decision variables can be expressed as physical quantities (e.g., production volumes), and nonlinear performance metrics such as Return on Investment (ROI) may be requested. These modeling choices introduce challenges, including the non-additivity of the objective function. This often results in highly concentrated optimized portfolios and thus limited diversification, which can be problematic for decision-makers seeking balanced investment strategies.
This paper proposes two strategies to enhance diversification in ROI-based PO models, both based on the Herfindahl–Hirschman Index (HHI). The first incorporates an HHI term directly into the objective function, with its corresponding weight allowing control over diversification. The second directly maximizes diversification while controlling expected profit and risk degradation around the optimum portfolio (obtained through conventional PO). Both strategies are evaluated using synthetic data (energy assets) to illustrate their behavior and practical trade-offs. The results highlight how each method can support different decision-making needs and enhance portfolio robustness.


\end{abstract}

\section{Introduction}

Portfolio Optimization (PO) is a key method for allocating investments across assets by balancing expected profit and risk \cite{daniel2025portfolio}. It typically represents a multi-objective optimization problem whose solution is a set of ‘non-dominated’ portfolios lying along an efficient frontier (or Pareto front).
In practice, PO is often implemented as a mono-objective optimization combining weighted expected profit and risk, and the efficient frontier is sampled by resolving the problem for various weight values \cite{daniel2025portfolio}.  
As a further step, decision-makers will choose their favorite portfolio within the set obtained based on additional criteria \cite{grodzevich2006normalization}.

In the finance literature, PO decision variables are typically expressed as percentages, such as financial asset shares, and the performance metric is often linear, measuring the return. However, in contexts involving 'physical assets', decision variables can be expressed in terms of physical quantities. For example, in an energy market application, the decision variables may represent the energy volumes to be produced by energy assets 
in which the decision-maker plans to invest \cite{rebennack2010energy}.
In addition, nonlinear performance metrics can be used, such as Return on Investment (ROI) \cite{mascomere2025unified}. These choices introduce challenges. First, they lead to a non-convex optimization problem (a question that we do not address here). Second, ROI is a non-additive metric and maximizing expected ROI often results in highly concentrated portfolios, strongly favoring the asset with the highest ROI. This behavior is undesirable, particularly because decision-makers seek greater diversification to mitigate the risk of concentration in a few assets. We found that adding a risk term to the objective function can improve diversification, but only to a limited extent. In contrast, when a linear performance metric is used, higher risk aversion is directly associated with greater diversification.
As a result, standard PO methods based on the ROI performance metric may fail to produce portfolios with the level of diversification expected by the decision-maker.


This paper addresses this limitation by proposing two strategies to enhance diversification in the context of PO with physical variables and non-linear performance metrics such as ROI. Both strategies are based on the Herfindahl–Hirschman Index (HHI) \cite{brezina2016herfindahl}. The first approach extends the objective function by adding an HHI term, whose weight provides control over diversification. The second directly maximizes diversification (i.e., minimizes HHI) while controlling expected profit and risk degradation around a given portfolio, which is obtained through conventional PO. Both strategies are evaluated using synthetic data (energy assets) to illustrate their behavior and practical trade-offs.

The paper is structured as follows: Section \ref{sec:Background} details the background and challenges related to our work. Section \ref{sec:Diversification} describes our proposed diversification strategies. Section \ref{sec:Results} presents results that highlight how each method can support different decision-making needs and enhance portfolio robustness. Section \ref{sec:Conclusion} concludes with key insights and future directions.

\section{Background}
\label{sec:Background}


\subsection{Mathematical Formulation}
\label{sec:Mathematical}

We consider single-period PO, which is commonly used in real-world applications. A portfolio is represented by a decision vector $\mathbf{x} = (x_1, \dots, x_n)$, where each component $x_i$ denotes the quantity of asset $i$ to purchase or produce given a budget $B$: 
\begin{equation}
\sum_{i=1}^n x_i = B
\quad,\quad
x_i\ge 0.
\label{eq:basic_const}
\end{equation}
In the case considered here, the $x_i$ can represent the amount to be produced across various energy assets, where $B$ denotes the total targeted production. Note that the variable $p_i=x_i/B$ is often used in the financial literature, leading to asset shares ($\sum_{i=1}^n p_i = 1$). Using $x_i$ with physical units is common for physical assets, but both formulations are equivalent.

The performance metric considered is ROI, calculated as:
\begin{equation}
f(\mathbf{x}) = \frac{\sum_{i=1}^n x_i \mathbf{R}_i}{\sum_{i=1}^n x_i \mathbf{I}_i},
\end{equation}
where $\mathbf{R}_i$ and $\mathbf{I}_i$ are the return and investment related to asset $i$ at the time of portfolio construction. In practice, these values are unknown. A statistical framework is adopted \cite{mascomere2025unified},  where $\mathbf{R}_i$ and $\mathbf{I}_i$ are treated as random variables (denoted in bold), and $m$ scenarios (indexed by $s$) are generated and assumed to be independent and identically distributed realizations of these random variables. This assumption is desirable because returns and investments are unknown at portfolio construction. Relying on a single realization would be misleading, so scenario distributions are used to estimate expected ROI and risk.

An estimate of the expected ROI (our profit measure) is then obtained using the $m$ scenario realizations:
\begin{equation}
 \text{ROI}(\mathbf{x}) = \frac{1}{m} \sum_{s=1}^m f^{(s)}(\mathbf{x})
\quad,\quad
f^{(s)}(\mathbf{x}) = \frac{\sum_{i=1}^n x_i R^{(s)}_i}{\sum_{i=1}^n x_i I^{(s)}_i}.
\label{eq: ROI}
\end{equation}

In terms of expected risk, we consider a measure that limits the loss-side variance. The semi-deviation proposed by Markowitz \cite{markovitz1959portfolio} relies on a Gaussian assumption. The Conditional Value-at-Risk (CVaR) deviation relaxes this assumption and captures loss-side tail features of the distribution more accurately, offering a more robust risk metric \cite{rockafellar2002deviation}:

\begin{align}\label{eq:CVaR-dev}
    \text{Risk}_\beta(\mathbf{x}) &= \text{ROI}(\mathbf{x}) - \min_\alpha F_\beta(\mathbf{x}, \alpha), \quad  \text{Risk}_\beta(\mathbf{x}) \ge 0\\
    F_\beta(\mathbf{x}, \alpha) &= \alpha + \frac{1}{1-\beta} \cdot \frac{1}{m} \sum_{s=1}^m \text{ReLU}(f^{(s)}(\mathbf{x}) - \alpha),\nonumber
\end{align}
where $\beta$ represents the confidence level and $\text{ReLU}(\cdot)$ denotes the rectified linear unit function.

In practice, PO is often formulated as a weighted trade-off between maximizing profit and minimizing risk:
\begin{equation}
    \underset{\mathbf{x} \in \mathcal{X}}{\text{max}} \left[ (1-w) \cdot  \text{ROI}(\mathbf{x}) - w \cdot \text{Risk}_\beta(\mathbf{x}) \right],
    \label{eq: obj function}
\end{equation}
where $\mathcal{X}$ represents the search space:
\begin{equation}\label{eq:constr}
    \mathcal{X}=\left\{ \mathbf{x} \Big| \sum_{i=1}^n x_i = B, x_i\ge 0, \text{other potential constraints}  \right\}.
\end{equation}
$\mathcal{X}$ can include potential operational constraints such as maximum production per asset, regulatory constraints, etc.

The set of non-dominated solutions along the efficient frontier is built by solving Eq. (\ref{eq: obj function}) for various values of the risk-aversion parameter $w \in [0,1]$.

\subsection{Challenges}

Two main challenges arise with ROI-based PO. The first is non-convexity, which we do not address here. The second, central to this article, is the tendency for the solution to concentrate all investment in the asset with the highest average ROI when no 'other potential constraints' are included in Eq. (\ref{eq:constr}) and when $w=0$ (i.e., only expected ROI) in Eq. (\ref{eq: obj function}). Indeed, we can show that:
%
%
\begin{equation*}
\mathrm{ROI}(\mathbf{x}) \le \max_j \left( \frac{1}{m} \sum_{s=1}^m \frac{R^{(s)}_j}{I^{(s)}_j} \right),
\end{equation*}
This implies that, in the absence of additional constraints and a risk measure, the maximum ROI solution is always fully concentrated in a single asset, i.e., leading to $p_i=\delta_{i,j^*}$ with $j^*=\arg\max_j ( \frac{1}{m} \sum_{s=1}^m R^{(s)}_j/I^{(s)}_j)$. Even when constraints and risk measures are added, ROI-based optimization often results in limited diversification. To address this issue, we propose two strategies that explicitly promote diversification in ROI-based PO models.

To the best of our knowledge, no existing work addresses the combined challenges of using ROI as a nonlinear performance metric, physical decision variables (e.g., production volumes), and their impact on portfolio diversification. This gap underscores the need for models that more accurately reflect the operational realities of the energy sector.

\section{Proposed Diversification approaches}
\label{sec:Diversification}

\subsection{Herfindahl–Hirschman Index (HHI)}


The HHI, used as the diversification metric, is defined as \cite{brezina2016herfindahl}:
\begin{equation}
\text{HHI}(\mathbf{x}) = \sum_{i=1}^{n} \left(\frac{x_i}{B}\right)^2,
\end{equation}
where \( \frac{x_i}{B} \) is the portfolio share of asset \( i \). The HHI ranges from \( \frac{1}{n} \) (maximum diversification, where all assets have equal portfolio shares) to 1 (fully concentrated portfolio).
Fostering diversification can thus be achieved by minimizing $\text{HHI}(\mathbf{x})$ or by maximizing $-\text{HHI}(\mathbf{x})$.

\subsection{Controlling Diversification via Objective Function Extension}
\label{subsec: Objective Function Extension}

The first proposition is to extend the PO model by including a HHI term:
\begin{equation}
    \underset{\mathbf{x} \in \mathcal{X}}{\text{max}} \left[ (1-w) \cdot \text{ROI}(\mathbf{x}) - w \cdot \text{Risk}_\beta(\mathbf{x}) - w_d \cdot \theta_1 \cdot \text{HHI}(\mathbf{x})\right],
\end{equation}
where a new weight \( w_d \cdot \theta_1 \) is defined to control the influence of diversification compared to the original problem in Eq. (\ref{eq: obj function}).
$w_d$ represents a diversification control parameter ranging from 0 to 1, where $w_d = 0$ ignores diversification and $w_d = 1$ corresponds to the maximum diversification allowed by the constraints (or very close to it). $\theta_1$ is introduced to rescale the HHI contribution to match the physical unit and order of magnitude of the ROI and risk terms. $\theta_1$ varies dynamically based on the value of $w$:
%
\begin{equation}
\theta_1 = \frac{w  \left|  \text{ROI}(\mathbf{x}) \right|_{mean}  + (1-w)  \left|  \text{Risk}_\beta(\mathbf{x}) \right|_{mean}}{ \left| \text{HHI}(\mathbf{x}) \right|_{mean}}.
\end{equation}

where $\left|\cdot\right|_{\text{mean}}$ denotes the average value of the corresponding metric (ROI, Risk, or HHI) computed over the set of non-dominated solutions along the Pareto front of the baseline optimal solution.

\subsection{Maximizing Diversification While Controlling Profit and Risk Degradation via Constraints Around an Optimum Portfolio}
\label{subsec: Controlling Degradation via Constraints}

The second strategy considers an objective function that focuses only on diversification, treating profit and risk as constraints, and allowing small deviations from the baseline solution $\mathbf{x}^*$ of Eq. (\ref{eq: obj function}). This enables control over acceptable degradation levels in each metric, aiming to obtain the most diversified portfolio within those limits.

The formulation is:
\begin{align}\label{eq:HHIcost}
    \underset{\mathbf{x} \in \mathcal{X}\cap\mathcal{X'}}{\text{min}} \text{HHI}(\mathbf{x}),
\end{align}
where $\mathcal{X}$ is defined by Eq. (\ref{eq:constr}) and:
\begin{eqnarray}\label{eq:constr2}
    \mathcal{X'}=\Big\{
    \mathbf{x} \Big|
    &&\text{ROI}(\mathbf{x}) \geq  \text{ROI}(\mathbf{x}^*)( 1- \Delta_p); \\
    &&\text{Risk}_\beta(\mathbf{x}) \leq \text{Risk}_\beta(\mathbf{x}^*) (1- \Delta_r)\Big\}.\nonumber
\end{eqnarray}
Unlike the previous approach, this formulation requires no weighting parameters. Instead, it introduces two new parameters, \( \Delta_p \) and \( \Delta_r \) (in $[0,1[$), which represent acceptable average profit and risk degradation percentages. When \( \Delta_p = \Delta_r = 0 \), the solution to Eq. (\ref{eq:HHIcost}) will most likely coincide with $\mathbf{x}^*$. In rare cases, another portfolio may be selected if a more diversified one (according to the HHI criterion) achieves exactly the same profit and risk values as $\mathbf{x}^*$, not just the same objective function value, but this is unlikely. 
This scheme can therefore be viewed as providing a controlled perturbation of the solution to Eq. (\ref{eq: obj function}).

An important point is that Eq. (\ref{eq:HHIcost}) is difficult to solve when the CVaR-deviation risk measure is used, whose efficient implementation involves a minimization over an auxiliary variable \( \alpha \) as in Eq. (\ref{eq:CVaR-dev}). In other words, the set in Eq. (\ref{eq:constr2}) actually is:

\begin{eqnarray}\label{eq:constr3}
    \mathcal{X'}=\Big\{
    \mathbf{x} \Big|
    &&\text{ROI}(\mathbf{x}) \geq  \text{ROI}(\mathbf{x}^*)( 1- \Delta_p);\\
    && \text{ROI}(\mathbf{x}) - F_\beta(\mathbf{x}, \alpha^*_\mathbf{x}) \leq \text{Risk}_\beta(\mathbf{x}^*) (1- \Delta_r);\nonumber\\
    &&\alpha^*_\mathbf{x}=\arg\min_\alpha F_\beta(\mathbf{x}, \alpha)
    \Big\}.\nonumber
\end{eqnarray}
As we search for `perturbations' of $\mathbf{x}^*$, the use of convex optimization tools such as Ipopt \cite{ipopt} is appropriate to solve these equations, as the local minimum closest to $\mathbf{x}^*$ would then be found, which is coherent. 
However, these equations remain impractical because the term $\arg\min_\alpha F_\beta(\mathbf{x}, \alpha)$ appearing in $\mathcal{X'}$, Eq. (\ref{eq:constr3}), cannot be directly put as a constraint in tools such as Ipopt  \cite{ipopt}.
Implementations involving interlaced iterations could be considered, but they would significantly increase the computational cost.

In the following, we propose an approximate formulation of Eqs. (\ref{eq:HHIcost})-(\ref{eq:constr3}) that do not lead to a strong overhead in the computational cost:

\begin{align}
    \underset{\mathbf{x} \in \mathcal{X}\cap\mathcal{X''}, \alpha}{\text{min}} \text{HHI}(\mathbf{x}) - w_r.\theta_2.
    F_\beta(\mathbf{x}, \alpha),
    \label{eq:fo_div_constr}
\end{align}
where
\begin{eqnarray}\label{eq:constr4}
    \mathcal{X''}=\Big\{
    \mathbf{x} \Big|
    &&\text{ROI}(\mathbf{x}) \geq  \text{ROI}(\mathbf{x}^*)( 1- \Delta_p);\\
    && \text{ROI}(\mathbf{x}) - F_\beta(\mathbf{x}, \alpha) \leq \text{Risk}_\beta(\mathbf{x}^*) (1- \Delta_r)
    \Big\}.\nonumber
\end{eqnarray}
\makebox[\linewidth]{This \hfill formulation \hfill  has \hfill the \hfill  advantage \hfill of \hfill allowing \hfill the \hfill computation \hfill of}
$\alpha^*_\mathbf{x} = \arg\min_\alpha F_\beta(\mathbf{x}, \alpha)$, which originally appears in $\mathcal{X'}$, Eq. (\ref{eq:constr3}), through an additional term in the objective function. 
This makes the scheme compatible with tools such as Ipopt but represents, in the general case, an approximation of the original scheme, as the HHI objective function is now modified by an additional term (the scheme will generally not converge to exactly the same decision vector). The approximation will be good if this modification is kept small, which is controlled by the weight \( w_r.\theta_2 \).
%
The \( \theta_2 \) term is introduced to align the scales of the HHI and $F_\beta$ terms, as illustrated below:
\begin{align}
    \theta_2 = \frac{\text{HHI}(\mathbf{x}^*)}{\text{ROI}(\mathbf{x}^*)-\text{Risk}_\beta(\mathbf{x}^*)},
    \label{eq:fo}
\end{align}
and $w_r$ is a very small, unit-less weight, tuned to keep the contribution of the $F_\beta$ term small compared to that of the HHI term.

\begin{figure}[h]
    \centering
    \includegraphics[width=0.6\textwidth]{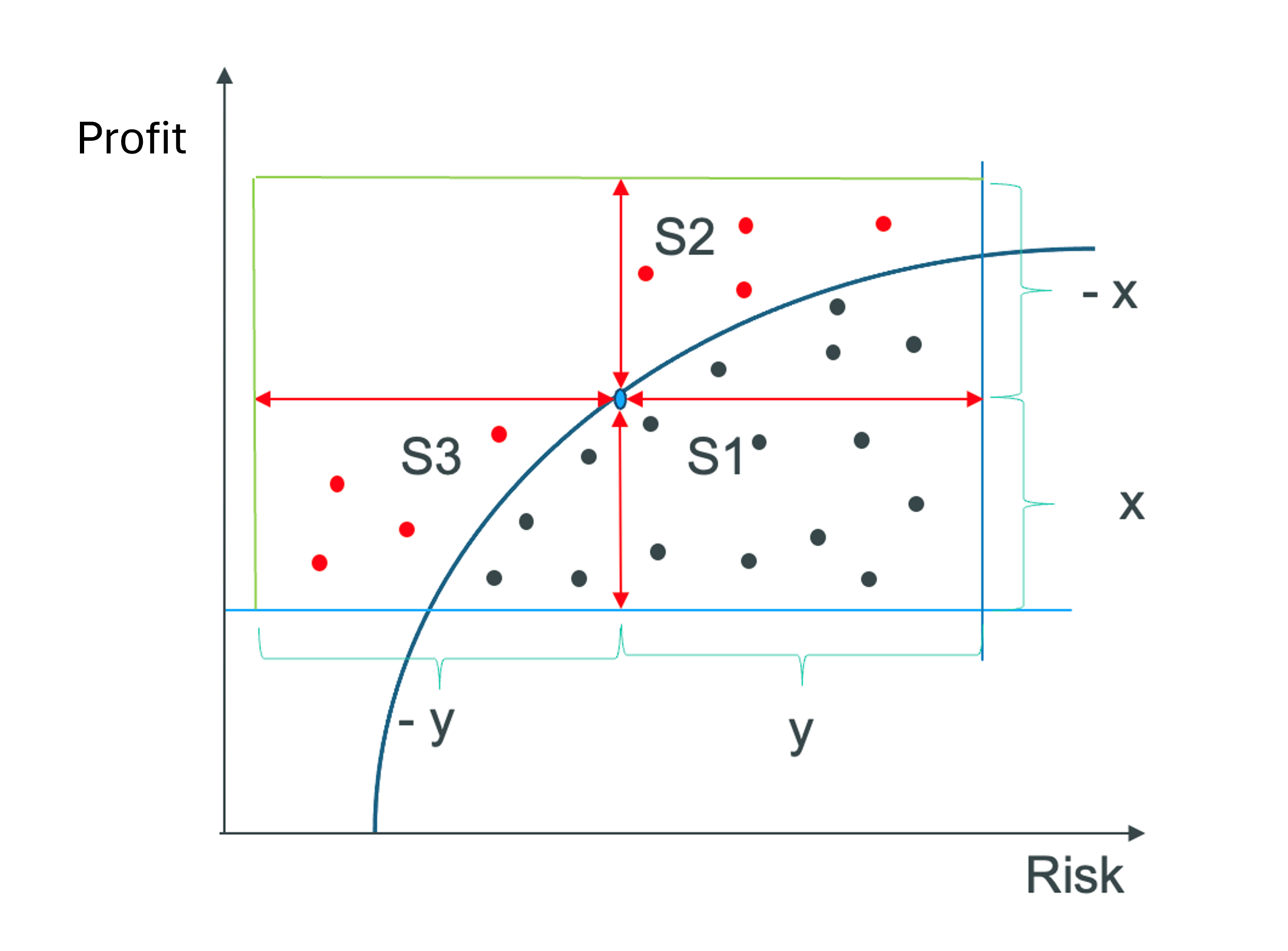}
    \caption{Illustration of the heuristic for generating tolerance pairs $(\Delta_p, \Delta_r)$ around an optimal portfolio $\mathbf{x}^*$ (blue point). The figure represents the efficient frontier (dark blue curve) in the Profit ($\mathrm{ROI}$ in our case) versus $\mathrm{Risk}$ (CVaR deviation in our case) plane, with $\mathbf{x}^*$ denoting the optimal solution obtained from the baseline PO, for a given $w$.}
    \label{fig:heuristic}
\end{figure}

Interestingly, this method enables the creation of multiple diversified solutions, or 'perturbed' portfolios compared to $\mathbf{x}^*$, by playing with various pairs of profit and risk tolerances $\Delta_p$, $\Delta_r$. Fig. ~\ref{fig:heuristic} illustrates the process. The rectangles $s_1$, $s_2$, and $s_3$ represent regions where random tolerance pairs are generated. Red points are infeasible, as they fall outside the efficient frontier.

The heuristic to generate a set of more diversified portfolio centered on the original PO result $\mathbf{x}^*$ follows these steps:
\begin{itemize}
    \item Run the baseline optimization using different values of $w$ to find optimal points in the Pareto front.
    \item For each selected point $(\mathrm{ROI(x^*)}, \text{Risk}_\beta(\mathbf{x}^*))$, define maximum acceptable degradation levels: $a$ for profit and $b$ for risk. 
    \item Construct a $2a \times 2b$ rectangle centered at the Pareto point. 
    \item Subdivide the area into four zones; three ($s_1$, $s_2$, $s_3$) are used for generating tolerance pairs.
    \item Specify the number of tolerance pairs to generate in each zone.
    \item In $s_1$, both profit and risk can degrade; in $s_2$, profit must improve and only risk may degrade; in $s_3$, risk must improve and only profit may degrade.
    \item Generate tolerance pairs $(\Delta_p, \Delta_r)$ randomly within each zone.
    \item Use each pair as input in the diversification-constraint model. These values define allowed degradation (or required improvement).
    \item Solve the model (Eqs. (\ref{eq:fo_div_constr})-(\ref{eq:constr4}) )for each tolerance pair, for each $w$ value used in the initial step.   
\end{itemize}

Some pairs, especially in $s_2$ and $s_3$, may yield infeasible solutions, as they can fall outside the efficient frontier. Feasibility can only be verified after solving each instance, since the frontier is not explicitly known a priori.

\section{Numerical Results}
\label{sec:Results}

\subsection{Energy Market PO}

PO in the energy sector has gained increasing attention due to the industry's structural complexity. This complexity stems from the integration of diverse generation technologies, varying regulatory frameworks, and market volatility. The literature distinguishes between short-term and long-term applications, each associated with specific modeling requirements and decision horizons.

Short-term PO models typically address operational decisions under uncertainty, such as bidding strategies in electricity markets \cite{narajewski2022optimal,faia2017ad}, real-time dispatch \cite{riddervold2021internal}, or short-term trading \cite{faia2015portfolio,faia2016optimization}. These models often incorporate stochastic prices, technical constraints, and policy considerations that affect market participation and profitability. 

In contrast, long-term PO studies are generally concerned with strategic investment planning \cite{tolis2011impact,selccuklu2023electricity}, and technology selection under demand growth, climate targets, and financial risk \cite{tolis2011impact,selccuklu2023electricity,maier2016risk}.

The following subsections present the application of the proposed diversification strategies to a synthetic case study involving long-term investment in energy-producing assets. The dataset is synthetic and 100 scenarios were used. Moreover, a simplified model was considered with only three constraints: maximum production per asset, maximum production per country, and a fixed CAPEX budget.

In the optimization model, each asset is defined as a unique combination of a technology and a country, and is assigned to either the Secured or Merchant category. Secured assets have low-variance scenario distributions, reflecting stable and predictable performance, while Merchant assets are characterized by broader distributions and greater market exposure. The results are presented using distinct colors for each asset, with labels following the format $T\{t\}\_C\{c\}\_\text{\{Category\}}$, where $T$ denotes the technology, $C$ the country, and $\text{Category}$ indicates either Secured or Merchant. $t$ and $c$ represent the indices of technology and country, respectively.


\subsection{Controlling Diversification via Objective Function Extension}

\begin{figure}[h]
  \centering
  \begin{subfigure}{0.49\textwidth}
    \includegraphics[width=\linewidth]{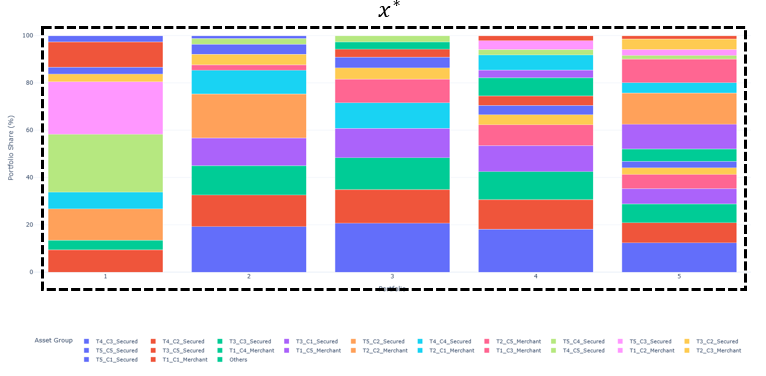}
    \caption{Optimization with $w_d = 0$}
    \label{fig:imagem1}
  \end{subfigure}
  \hfill
  \begin{subfigure}{0.49\textwidth}
    \includegraphics[width=\linewidth]{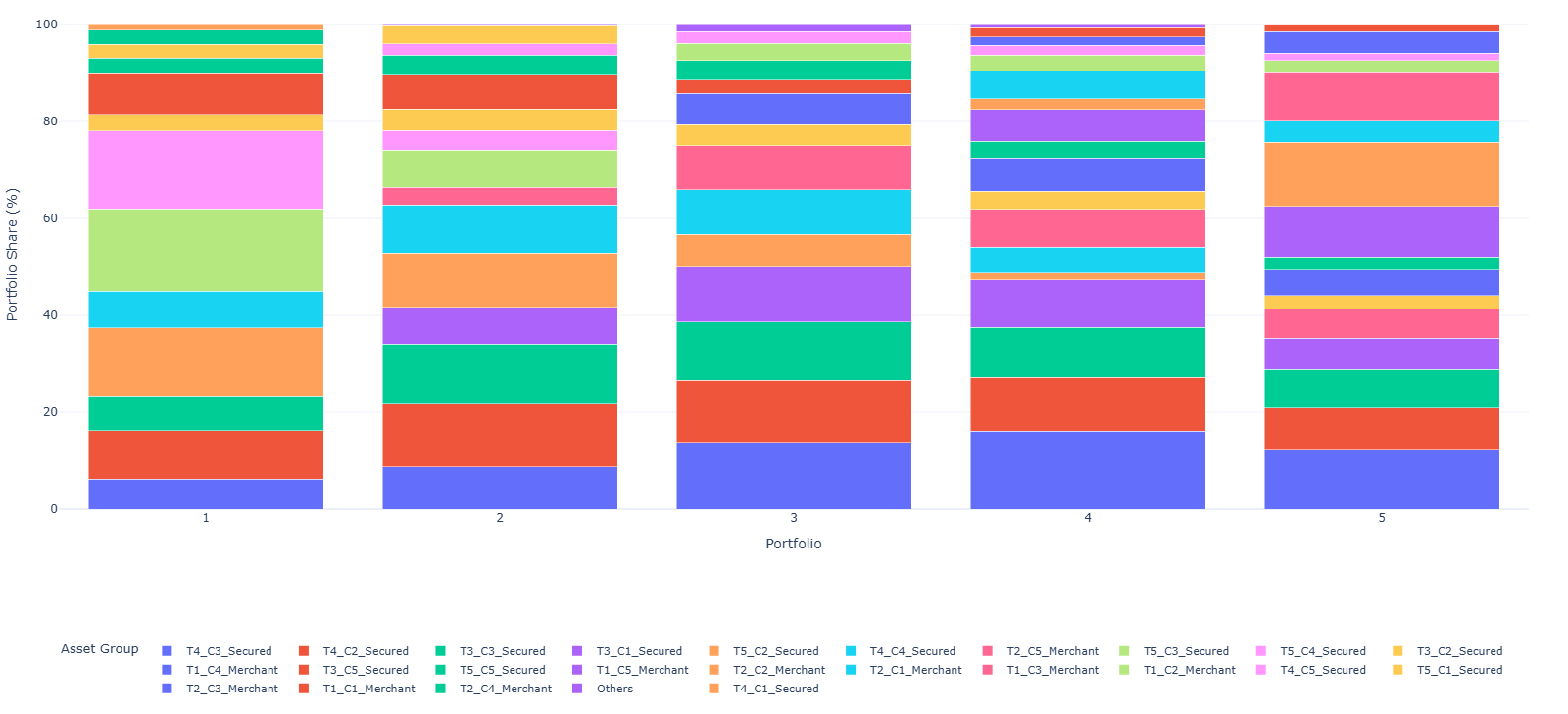}
    \caption{Optimization with $w_d = 0.2$}
    \label{fig:imagem2}
  \end{subfigure}
  \begin{subfigure}{0.49\textwidth}
    \includegraphics[width=\linewidth]{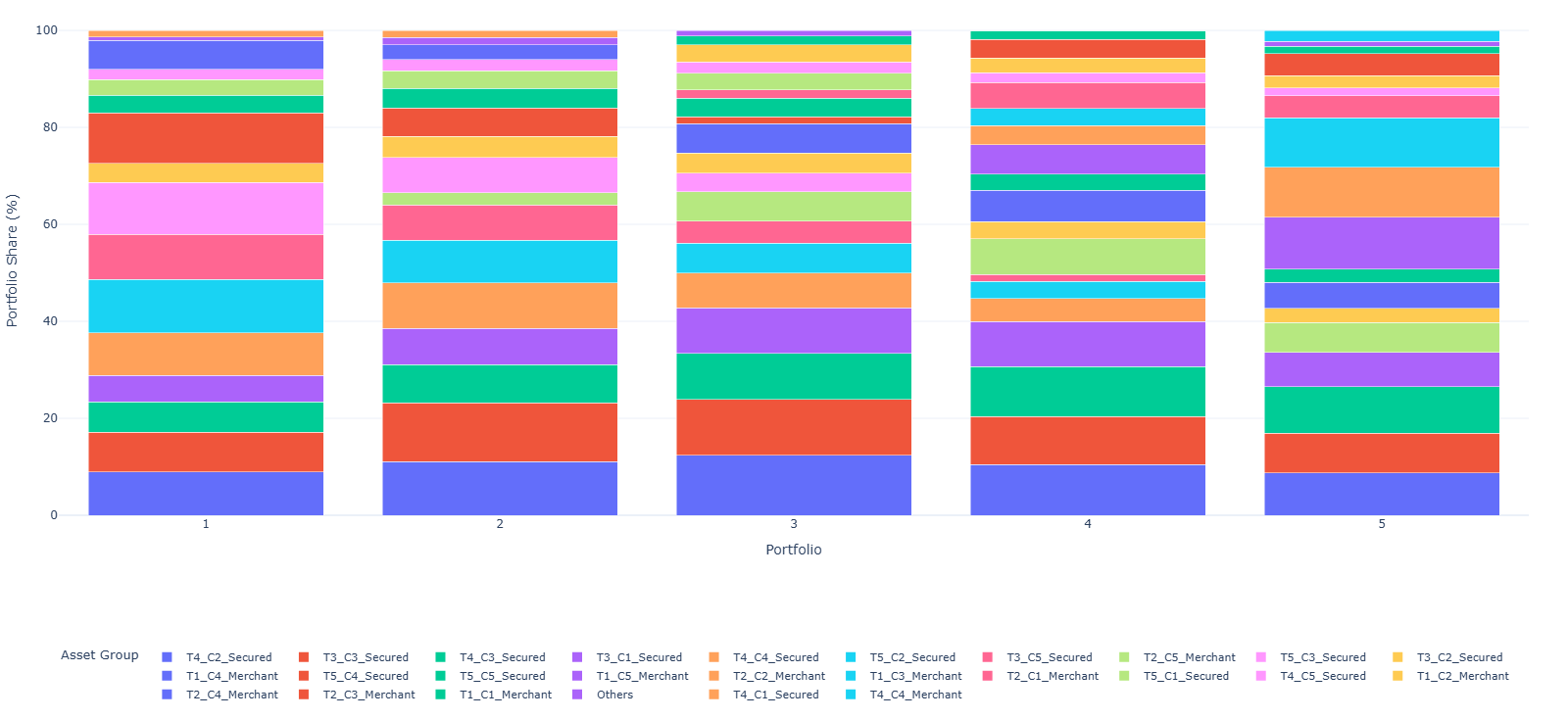}
    \caption{Optimization with $w_d = 0.5$}
    \label{fig:imagem3}
  \end{subfigure}
  \hfill
  \begin{subfigure}{0.49\textwidth}
    \includegraphics[width=\linewidth]{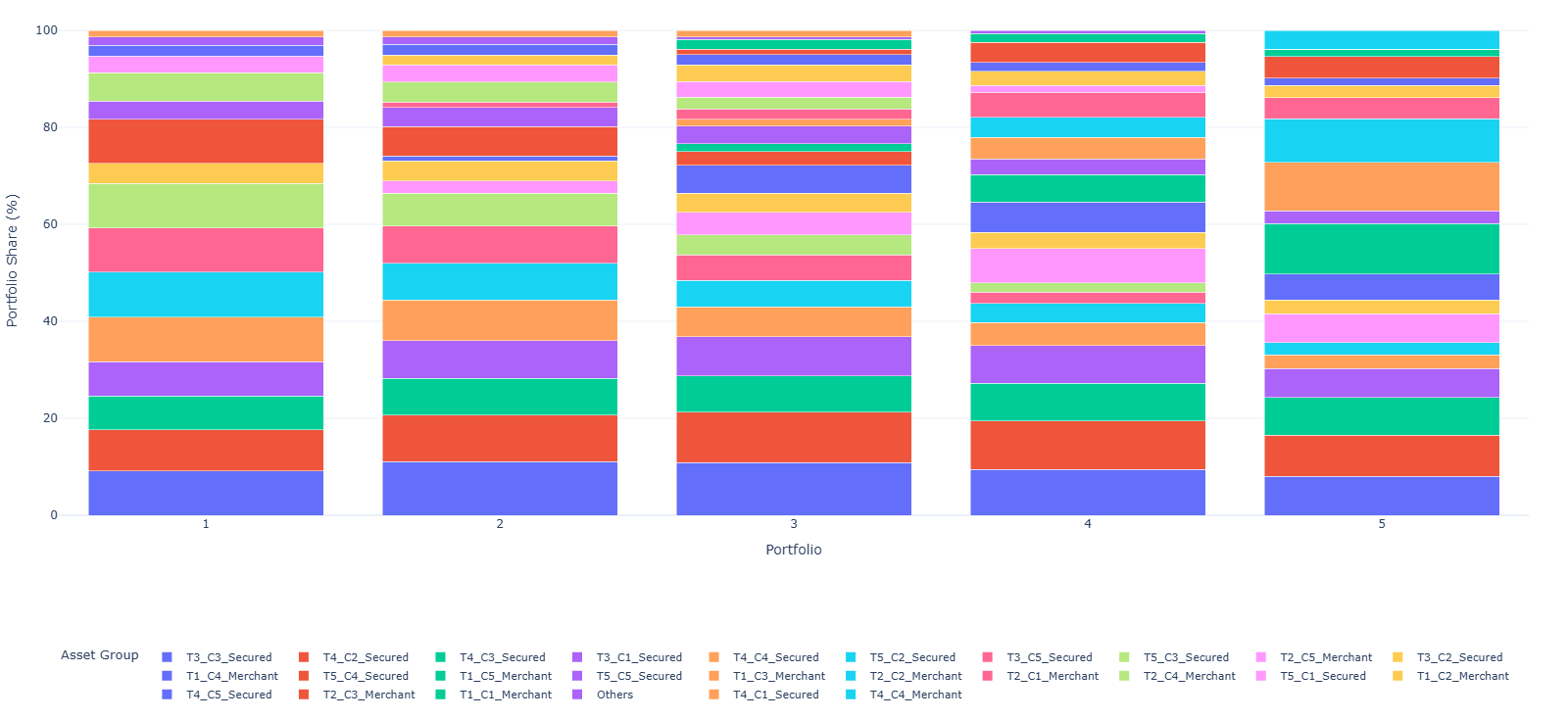}
    \caption{Optimization with $w_d = 0.9$}
    \label{fig:imagem4}
  \end{subfigure}
  \caption{Results of diversification via objective function extension. Panel (a) shows the baseline optimal portfolio $\mathbf{x}^*$ ($w_d = 0$). Panels (b)-(d) illustrate the effect of increasing $w_d$ (diversification weight). Each asset is represented by a distinct color.}
  \label{fig:div_obj_roi}
\end{figure}

To evaluate the effect of directly promoting diversification during optimization, we apply the strategy presented in subsection \ref{subsec: Objective Function Extension}. Fig.~\ref{fig:div_obj_roi} shows the results for this approach. Each sub-figure contains five columns representing the portfolio composition for five different values of \( w \): 1, 0.8, 0.6, 0.4, and 0.2. In Fig. \ref{fig:div_obj_roi}(a), diversification is not considered in the objective (\( w_d = 0 \)). Fig. \ref{fig:div_obj_roi}(b), \ref{fig:div_obj_roi}(c), and \ref{fig:div_obj_roi}(d) illustrate increasing diversification weights ($w_d = 0.2, 0.5, 0.9$), leading to progressively more diversified portfolios. This behavior is smooth and consistent across all $w$ values.

\begin{figure}[h]
    \centering
    \includegraphics[width=0.735\linewidth]{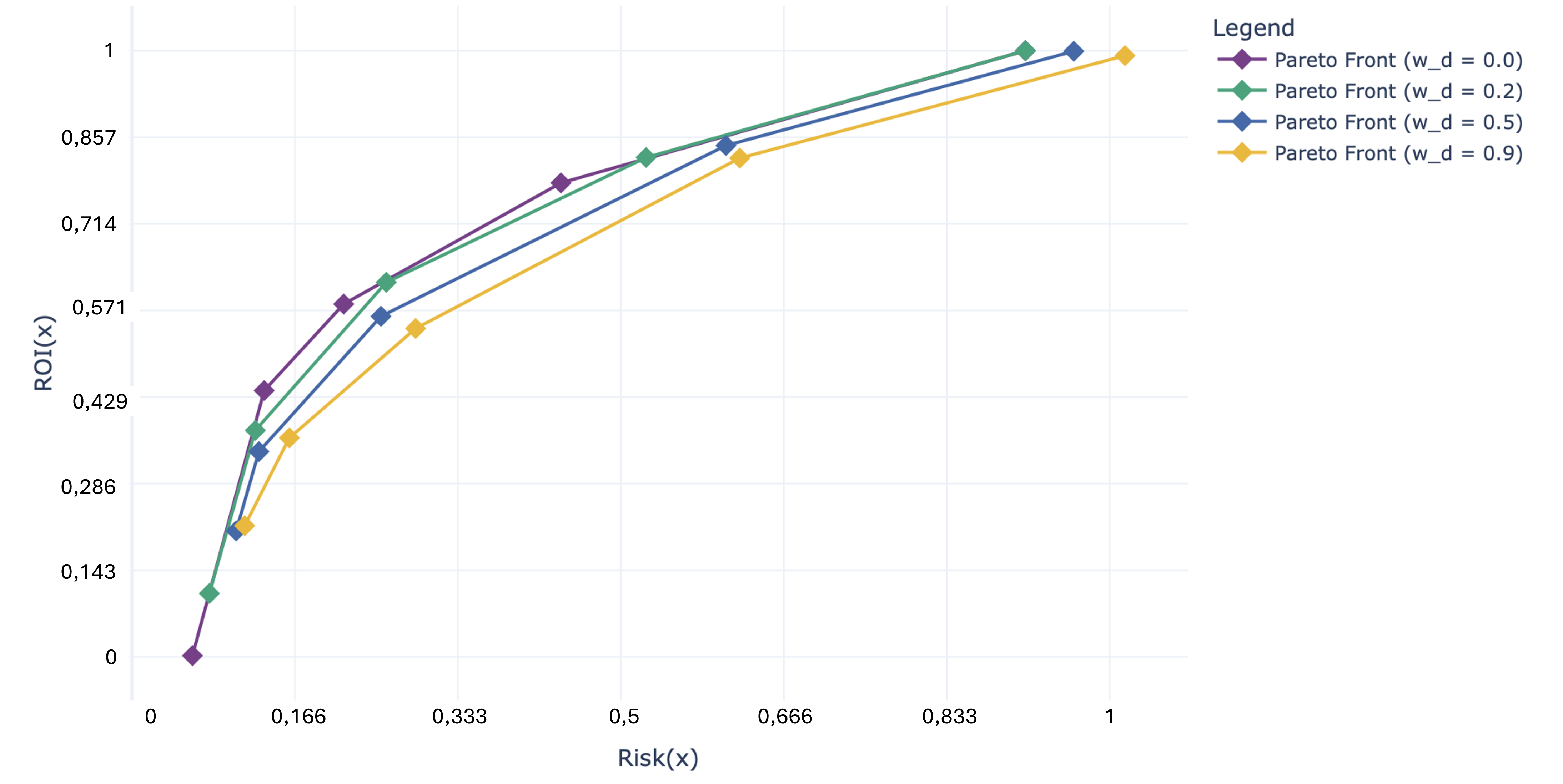}
    \caption{Comparison of efficient frontiers for different values of the diversification weight $w_d$ in the $\mathrm{ROI(\mathbf{x})}$ versus $\mathrm{Risk_\beta(\mathbf{x})}$ plane. Each curve corresponds to a specific $w_d$, with points representing solutions for various $w$. The axes were normalized to the interval [0, 1].}
    \label{fig:div_obj_roi_impact}
\end{figure}

However, greater diversification directly impacts ROI and risk, which can be evaluated a posteriori. Fig.~\ref{fig:div_obj_roi_impact} shows this trade-off. The Pareto front is plotted for different values of the diversification weight $w_d$:  purple ($w_d = 0$), green ($w_d = 0.2$), blue ($w_d = 0.5$), and yellow ($w_d = 0.9$), with five points per curve corresponding to decreasing $w$ values (from 1 to 0.2). 

At $w = 1$, both risk and ROI increase. For intermediate values of $w$ ($w = 0.8$ and $w = 0.6$), risk generally continues to rise while profit tends to decrease. For the lowest values of $w$ ($w = 0.4$ and $w = 0.2$), profit increases slightly. This behavior reflects the complexity introduced by the use of ROI as indicator. Unlike linear metrics, the relationship between numerator and denominator is neither direct nor proportional, making the effect of diversification on the objective function harder to interpret.
In contrast to linear models, where increased diversification typically leads to a natural reduction in risk, this relationship does not necessarily hold for the ROI formulation. 


\subsection{Maximizing Diversification while Controlling Profit and Risk Degradation via Constraints}

This subsection presents the results obtained by applying the strategy described in subsection \ref{subsec: Controlling Degradation via Constraints}. The optimal profit and risk values obtained from this baseline serve as references to define the degradation bounds in the diversification step.

Here, the same set of randomly generated tolerances was used for all values of $w$, selected within a maximum degradation of 10\% (i.e., in Fig. \ref{fig:heuristic}, $x = y = 0.1$). The weight in Eq. (\ref{eq:fo_div_constr}) was set to $w_r = 0.001$.

\begin{figure}[h]
    \centering
    \includegraphics[width=0.735\textwidth]{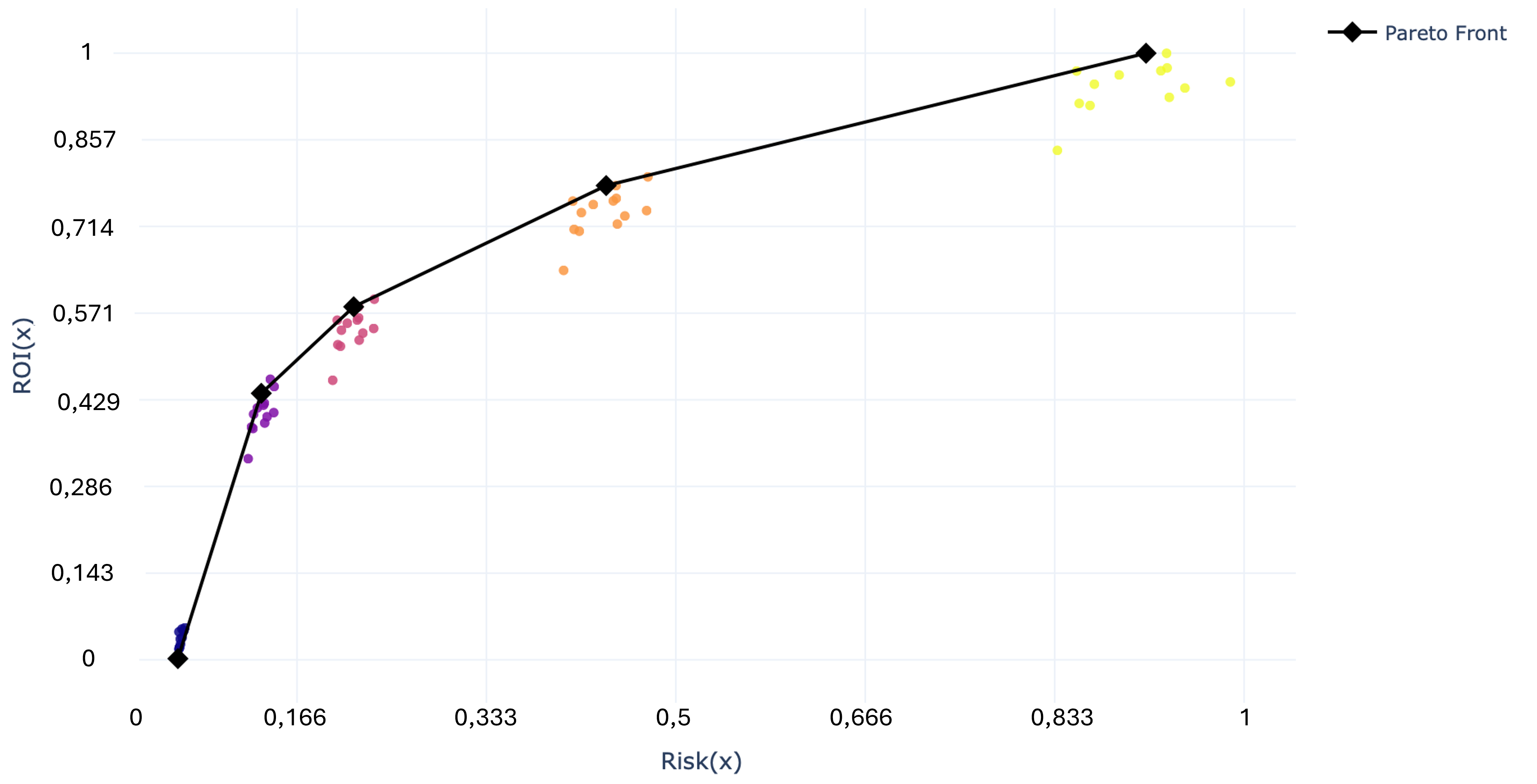}
    \caption{Optimal and diversified suboptimal points in the $\mathrm{ROI(\mathbf{x})}$ versus $\mathrm{Risk_\beta(\mathbf{x})}$ plane. The black points indicate the optimal portfolios on the efficient frontier for each value of $w$, while the colored points represent diversified suboptimal portfolios generated for different tolerance pairs $(\Delta_p, \Delta_r)$. The axes were normalized to the interval [0, 1].}
    \label{fig:heuristic_roi}
\end{figure}

Fig.~\ref{fig:heuristic_roi} displays the optimal portfolios on the Pareto front for each $w$, together with the diversified suboptimal portfolios obtained for each degradation pair $(\Delta_p,\Delta_r)$. The diversified solutions are similarly distributed around all optimal points, although with different proportions. This occurs because maximizing diversification tends to fully, or nearly fully, saturate the ROI and risk constraints, resulting in portfolios whose ROI and risk values are very close to the imposed degradation bounds.

\begin{figure}[h]
  \centering
  \begin{subfigure}{0.49\textwidth}
    \includegraphics[width=\linewidth]{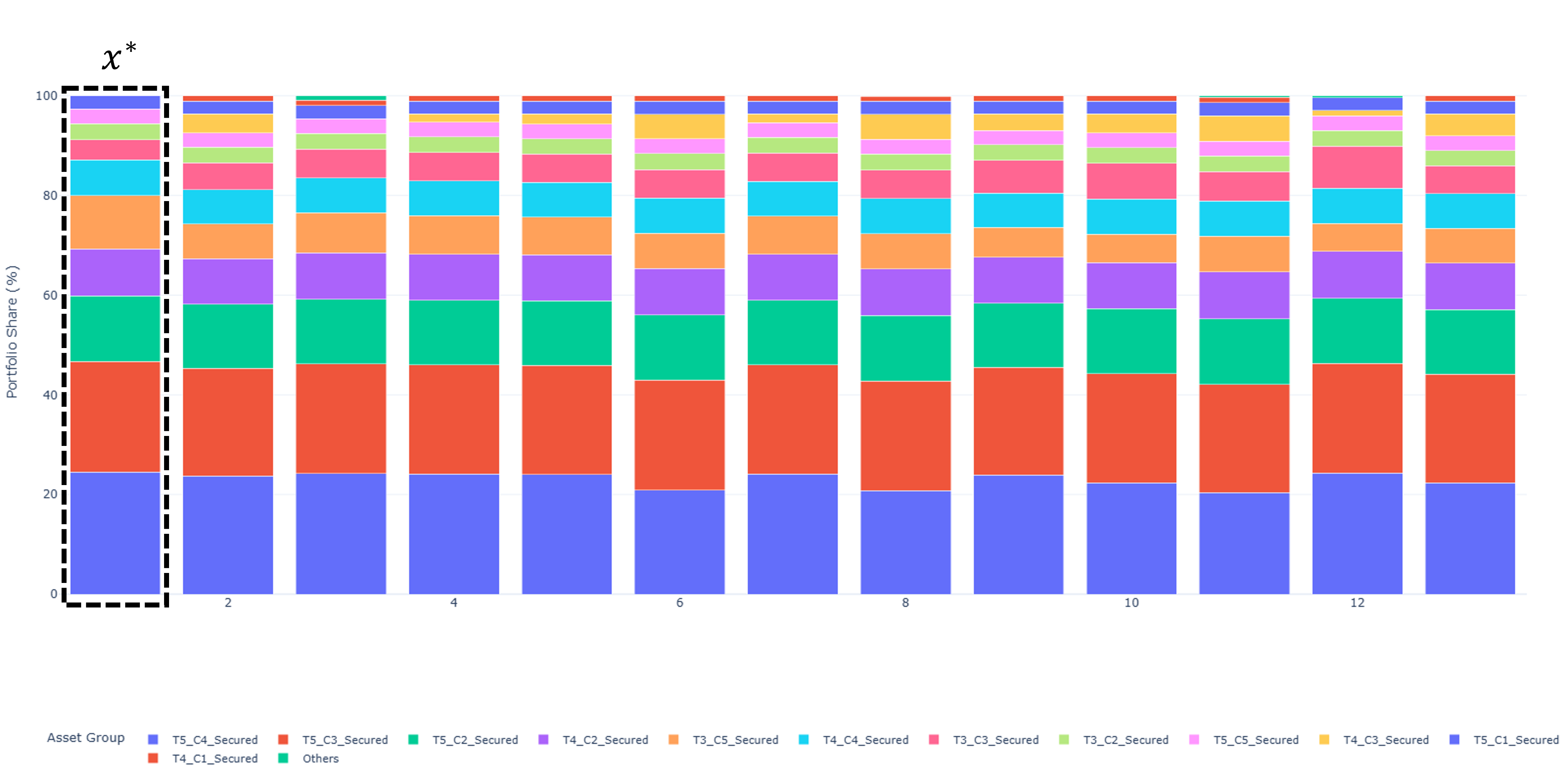}
    \caption{$w=1$}
    \label{fig:imagem6}
  \end{subfigure}
  \hfill
  \begin{subfigure}{0.49\textwidth}
    \includegraphics[width=\linewidth]{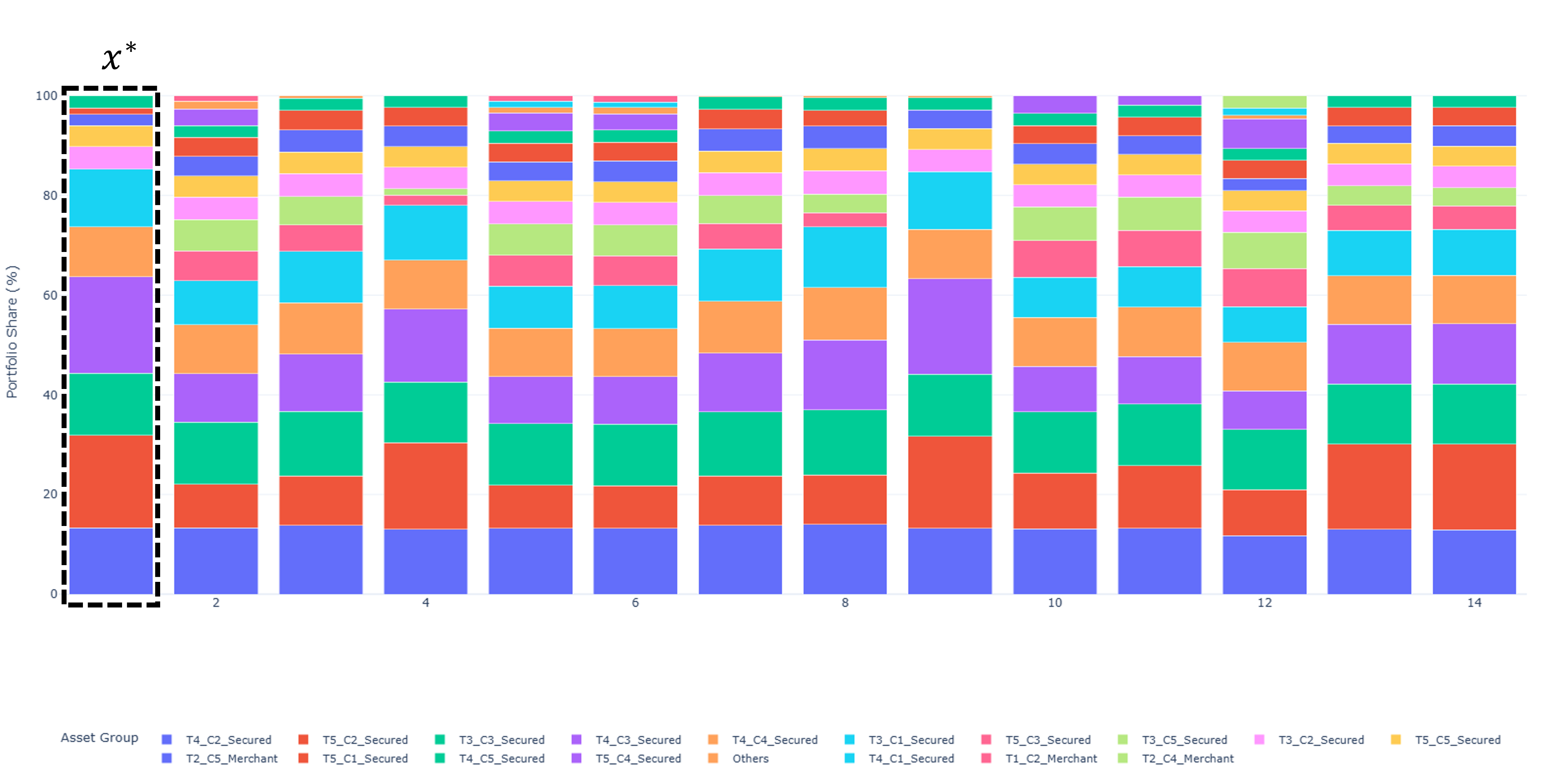}
    \caption{$w=0.8$}
    \label{fig:imagem7}
  \end{subfigure}
  \begin{subfigure}{0.49\textwidth}
    \includegraphics[width=\linewidth]{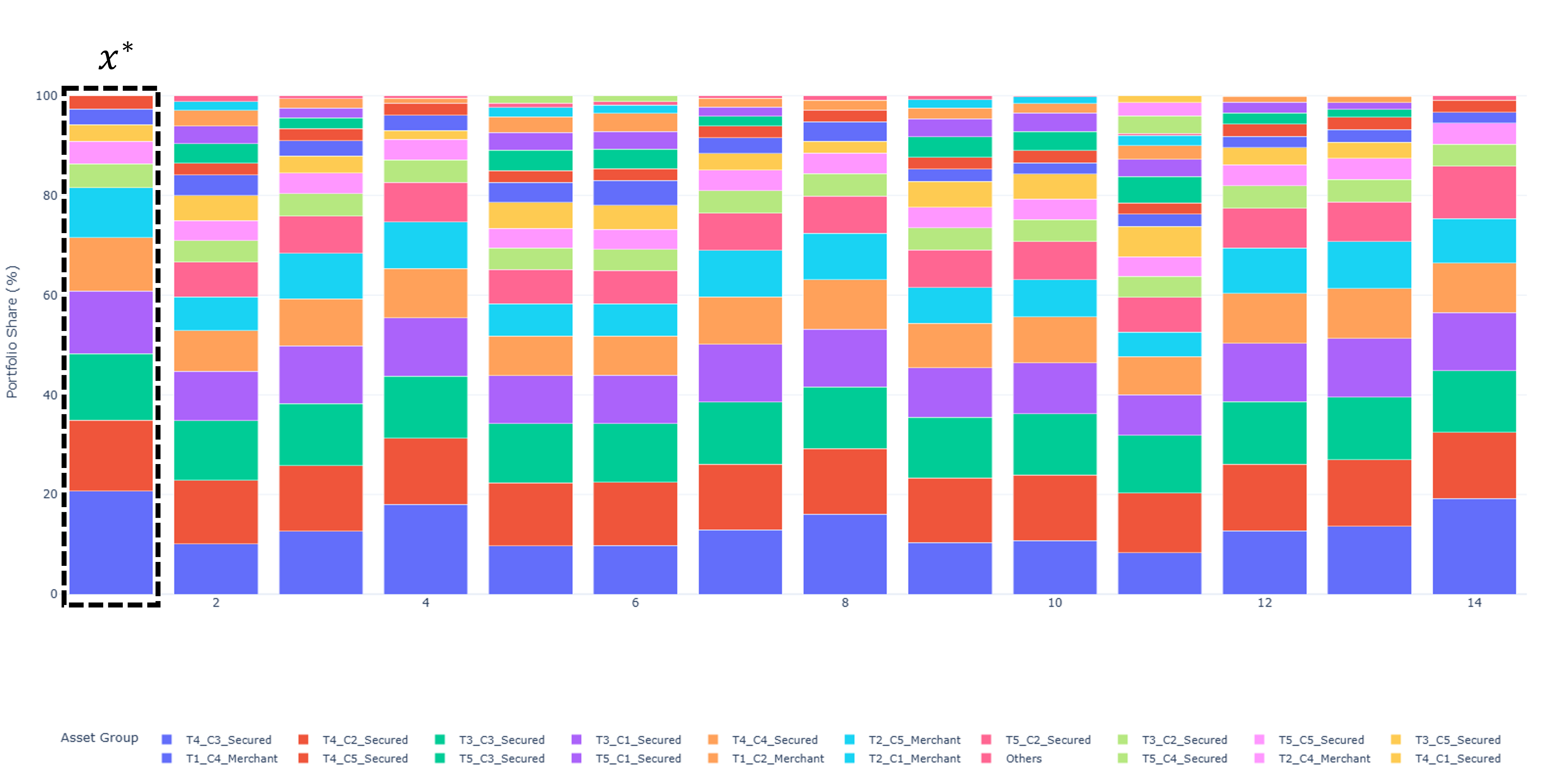}
    \caption{$w=0.6$}
    \label{fig:imagem8}
  \end{subfigure}
  \hfill
  \begin{subfigure}{0.49\textwidth}
    \includegraphics[width=\linewidth]{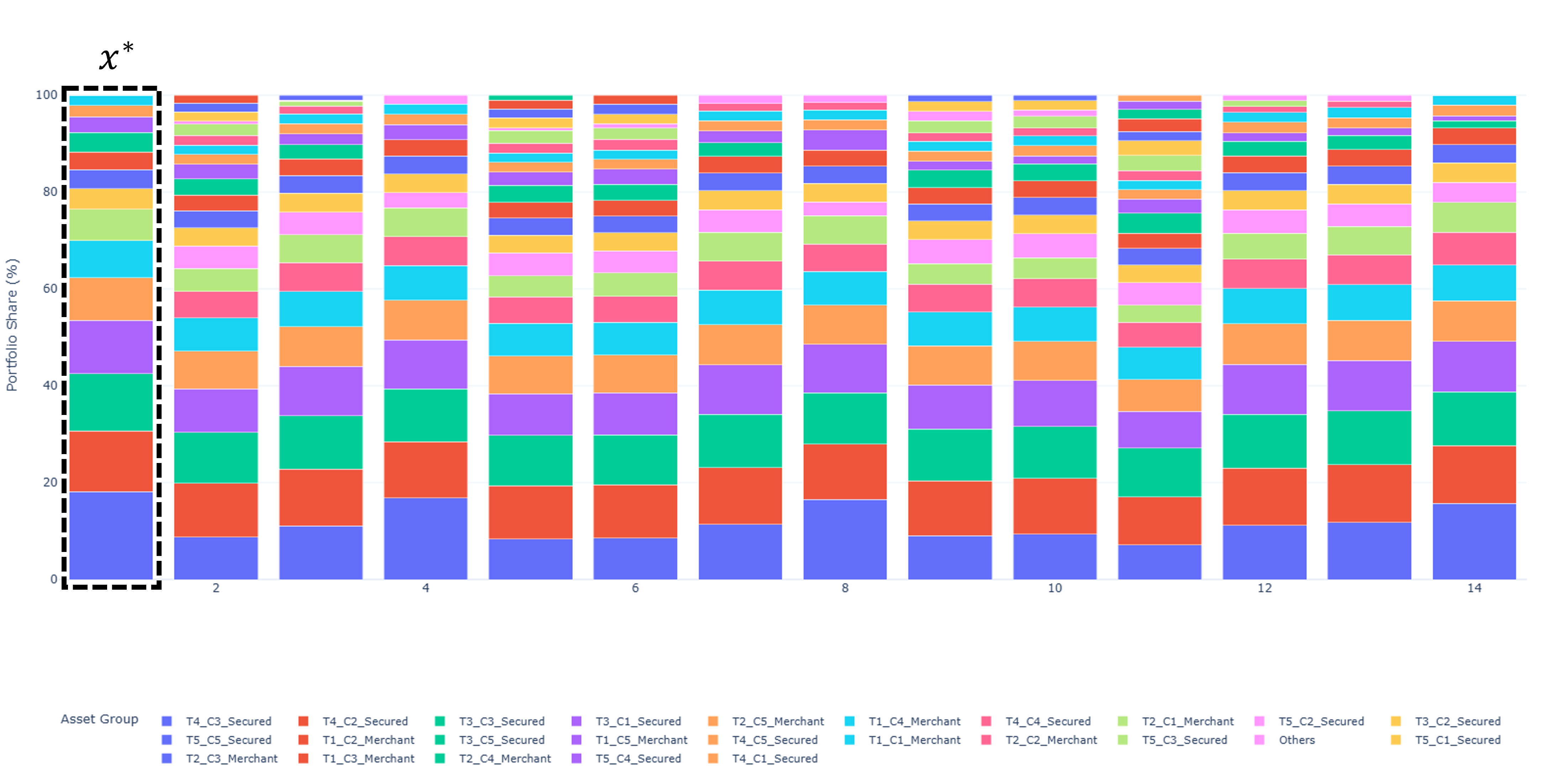}
    \caption{$w=0.4$}
    \label{fig:imagem9}
  \end{subfigure}
  \begin{subfigure}{0.49\textwidth}
    \includegraphics[width=\linewidth]{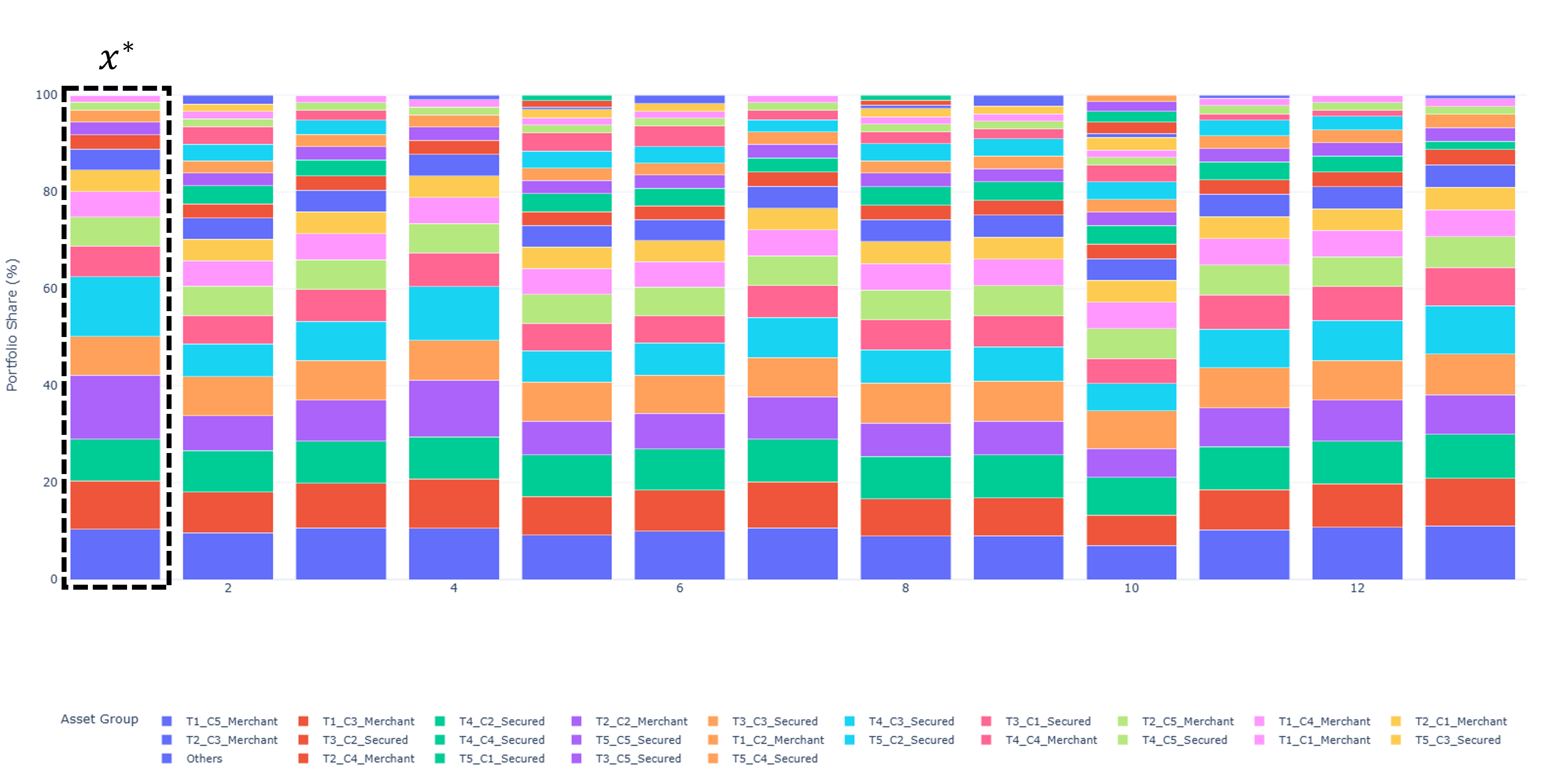}
    \caption{$w=0.2$}
    \label{fig:imagem10}
  \end{subfigure}
  \caption{Portfolio share distributions for optimal and suboptimal points. In each subfigure (corresponding to a different $w$), the optimal portfolio $\mathbf{x}^*$ is highlighted with a black dashed outline. Each asset is represented by a distinct color.}
  \label{fig:roi_portfolios}
\end{figure}

Fig.~\ref{fig:roi_portfolios} shows the portfolio compositions for each $w$, comparing the optimal portfolios with the corresponding perturbed ones. In each sub-figure, the first column corresponds to the optimal solution, while the remaining columns represent diversified portfolios generated by the heuristic.

Finally, we observe that some assets remain constant or nearly unchanged across all solutions, even when diversification is explicitly encouraged. These assets appear to be essential for maintaining portfolio profitability or controlling overall risk.

\section{Conclusions and Perspectives}
\label{sec:Conclusion}

We investigated PO models that use physical decision variables and nonlinear performance metrics, such as ROI, with an application in the energy sector. Two diversification strategies based on HHI were proposed and demonstrated using synthetic data.
The first strategy incorporates a HHI term directly into the objective function, allowing the decision-maker to control the desired level of diversification through a weight. This method offers a direct way to explore trade-offs between profit, risk, and diversification. The second strategy maximizes diversification (minimizes HHI) under constraints that limit degradation in profit and risk relative to an optimal portfolio obtained by conventional PO. It offers the possibility to generate multiple diversified portfolio perturbations, providing additional options to the decision-maker.
Both strategies offer complementary advantages, and the choice depends on the decision-maker's priorities. 

Future work may focus on parallelizing the optimization process and using the resulting portfolio sets for confidence interval analyses. We also plan to apply the method presented in Section 4.4 of \cite{mascomere2025unified} to enable a homogeneous comparison of PO feature variations (here, the portfolio perturbations obtained with our proposed methods) through marginal costs.

\begin{credits}
\subsubsection{\ackname} We thank Matthieu Blondel, Thomas Balsan, and Emeric de Monteville for their valuable feedback. We also acknowledge Jean-Patrick Mascomère and Yagnik Chatterjee for their support during the project. Finally, we thank TotalEnergies for permission to publish this work.

\end{credits}
%
%
%
\bibliographystyle{splncs04}
\bibliography{mybibliography}
%

\end{document}